\input amstex
\documentstyle{amsppt}
\headline={\hss\tenrm\hss\folio}
\vsize=9.0truein
\hsize=6.5truein
\parindent=20pt
\baselineskip=22pt
\NoBlackBoxes

 \topmatter
 \title Orbifold Index Cobordism Invariance \endtitle 
 \author Carla Farsi, \
 Department of Mathematics, \ 
 University of Colorado, \
 395 UCB, \
 Boulder, CO 80309--0395, USA. \ 
 e-mail: farsi\@euclid.colorado.edu 
 \endauthor
 \abstract 
We prove cobordism index invariance for 
pseudo-differential elliptic 
operators on closed orbifolds
with $K$--theoretical methods.
\endabstract
\endtopmatter

 \noindent {\bf 0.  Introduction.}

Orbifolds, which   play 
an important role 
in mathematical physics, have been recently 
studied from many viewpoints,
including groupoid and symplectic geometrical. 
Orbifold index theory has played an 
important role in these studies.
 
Firstly, Atiyah  established in [At] many 
important results on the index 
of $G$--transversally elliptic operators, 
which are closely related to elliptic operators on 
orbifolds as we will explain later.
In the late seventies   
 Kawasaki  gave several  
 proofs,  one of which making extensive 
 use of Atiyah's results, of an index theorem for
 orbifolds, [Kw1], [Kw2], [Kw3]. I  
 later proved a $K$--theoretical index theorem
 for orbifolds with operator algebraic means, [Fa1], 
 and  Berline and Vergne computed the 
 index of transversally elliptic
 operators via heat equation methods in [BV]. 
 Their work was then deepened by Vergne who proved 
 a  general index theorem
 for orbifolds, [V].  By using elliptic estimates,  I
 established  in [Fa2] some spectral
 properties of the  eigenvalues of the Laplacian on orbifolds, 
 and in [Fa3] 
 I defined orbifold eta invariants and established an   
 index theorem for orbifolds with boundary.

There have been many 
proofs of the cobordism invariance of the 
index of pseudo-differential elliptic operators on closed manifolds.
(See for example [Pa], [Hi], [Br1].)  Recently, Carvalho proved 
in [Ca] a $K_G$--theoretical ($G$ compact Lie) index cobordism 
invariance theorem for pseudo--differential 
$G$--equivariant elliptic operators
that are multiplication at infinity. 
Her methods are topological and rely heavily on key 
properties of Atiyah's $K_G$--functor.
It has also recently come to our attention that
Braveman proved the  cobordism invariance of
orbifold indices analytically, see [Br2].

In this note we will prove the cobordism invariance of the index of 
pseudodifferential operators on orbifolds topologically by a
generalization of Carvalho's method. Our main result is indeed the following theorem.

\proclaim{Theorem  4.6} Let $Q_i$, be a closed 
 orbifold which is  the  locally 
free quotient of an action of a compact Lie group
on a smooth manifold, and   let $P_i$ be an elliptic pseudo-differential operator on $Q_i$ with symbol $p_i$, 
$i=1,2$. Suppose that $(Q_1, p_1)$ is orbifold symbol cobordant to $(Q_2, p_2)$. Then  $Ind(P_1)=Ind (P_2)$.
\endproclaim

To prove this result, instead than working directly on the orbifold, we work on its $G$-frame bundle. 
More in general, any closed effective orbifold $Q$ 
arises as the locally free  quotient
of an action of a compact Lie group $G$ on a closed $G$--manifold
$M$, see e.g. [Fa1]. ($M$ is a $G$--frame bundle of $Q$, and in general $G$ is a compact Lie group.)
Therefore an elliptic pseudo--differential operator $P$  
on $Q$ lifts to a $G$--transversally 
elliptic operator ${\tilde P}$ on $M$  with $G$--transverse symbol class 
in $K_G^0(T^*_G(M))$. Atiyah's distributional index homomorphism $Ind:
K_G^0(T_G(M))\to \Cal D'(G)$ calculates the index of $P$, and, consequently, 
the index of $\tilde P$; see Section 1 for details.

 There are two main ingredients in our proof.  One is the push-forward property of the transverse index $Ind$ under $G$-embeddings established by Atiyah in [At], and the other is the $K$--theoretical proof provided by Carvalho in [Ca] for the invariance of the index of elliptic operators on manifolds under cobordism. Here we generalize Carvalho's approach to the context of transversally elliptic operators in the framework of [At]. So $Ind (\sigma)=0$ if $\sigma$ arises from a trivial symbol $G$--cobordism, see Theorem 4.3. 
By reinterpreting this at the orbifold level, 
we obtain our orbifold index cobordism invariance 
result, Theorem 4.6. 

More in detail, the contents of this note are as follows.  In Section 1
we recall the definition of Atiyah's distributional transverse index $Ind$. 
In Section 2, we define restriction and boundary maps,  
together with  symbol $G$--cobordism.
In Section 3 we detail properties of the index, restriction, and
boundary maps
with respect to $G$--equivariant embeddings. 
In Section 4 we will state and prove our main results.

 In the sequel, all orbifolds and manifolds are assumed to be smooth,
 $Spin^c$, connected, and closed, unless otherwise specified. Moreover, $G$ 
 will denote a compact connected  
Lie group. By a $G$--manifold we mean a manifold
 with a smooth and proper action of $G$. 

 \vskip 1em
 \noindent {\bf 1.  The Distributional Index Homomorphism.}
 \vskip 1em

 In this section we will review the definition and some of the main
 properties of 
 the distributional $G$--index homomorphism for $G$--transversally 
 elliptic operators  (c.f. [At] for details).

 Let $X$ be a $G$--manifold and let $\Cal G$ be the Lie algebra of $G$.
 To each $V\in \Cal G$ associate the vector field $V_G$ on $X$ defined
 by
 $$
 V_G(f)= {\hbox {Lim}}_{t\to 0} \frac{f({Exp}(tV_x))
 -f(x) }t, \quad \forall f\in C^\infty(X).
 $$

 \proclaim{Definition 1.1 [At]} Let $X$ be a $G$--manifold. Define the
 $G$--invariant space $T^*_G(X)\subseteq T^*(X)$ by
 $$
 T^*_G(X)=\left\{  v\in T^*(X)\quad | \quad v(V_G)=0\quad \forall\, V\in
 \Cal G  \right\}.
 $$
 \endproclaim

Let $X$ be a $G$--manifold and let $D$ be a pseudo--differential operator acting on sections of the $G$--vector bundle $E$. $D$ is said to be $G$--transversally elliptic $G$ if the symbol of $D$ is invertible on $T^*_G(X)$, except for the zero section. We will call such an operator a $G$--t.e.p.d. operator for short.

We will now recall how the transverse index of $D$ is defined, see [At; Lecture 2]. If $\Cal G$ denotes the Lie algebra of $G$, and $\Cal X_j$, $j=1,\dots, k$ the first order differential operators defined by the action of $\Cal G$ on $E$, denote by $\Delta_G $
the following operator  
$$
\Delta_G = 1-\Sigma_{j=1}^k {\Cal X_j^2}.
$$
Let $\lambda$ be an eigenvalue of $\Delta_G$, and 
denote by ${\Cal C}^\infty (X, E)_\lambda$ the kernel of 
the operator $\Delta_G -\lambda$. Since $D$ is $G$--invariant, $D$ commutes with
$\Delta_G$, and induces an operator
$$
D_\lambda : {\Cal C}^\infty (X, E)_\lambda \to {\Cal C}^\infty (X, E)_\lambda,
$$
with index $ Ind(D_\lambda) $. Define
$$
Ind(D) =\Sigma_\lambda Ind(D_\lambda).
$$
This sum converges in the sense of distributions, and is equal to  
the distibutional index of $D$. If we denote by $\Cal D'(G)$  the group of the $G$--invariant distributions on $G$,  then the index of $D$  is an element of  $\Cal D'(G)$, [At].

Let $K_G^s$, $s=0,1$, be Atiyah's equivariant $K$--theory functor. Then  
the symbol of $D$ determines a class $\sigma_D\in K_G^0(T^*_G(X))$, in analogy $G$--elliptic operators, [At].  $\sigma_D\in K_G^0(T^*_G(X))$ is called the $G$--transverse symbol class of $D$.  We have, 

 \proclaim{Theorem 1.2 [At; Theorem 2.6]}  Let $X$ be a $G$--manifold and let $D$ be a $G$--t.e.p.d. operator.  Then the index of $D$ depends only on the class
$\sigma_D\in K_G^0(T^*_G(X))$. In particular, there exists a index homomorphism
 $$
 Ind: \, K_G^0(T^*_G(X))\to \Cal D'(G),
 $$
 such that $Ind (\sigma_D)=Ind(D)$. 
 \endproclaim

 \noindent $Ind$ can also be defined for non--compact $G$--manifolds 
 via equivariant $G$--embeddings into  compacts $G$--manifolds, [At]. 

If the action of $G$ on $X$ is locally free, and if $D$ is the 
lift of a pseudo--differential  elliptic operator on the quotient
orbifold $X/G$, $Ind$ computes the distributional orbifold index, and consequently the orbifold numerical index of the operator  [Kw3], [Fa1], [V].

 \vskip 1em
 \noindent {\bf 2.  Boundary Maps and Symbol $G$--Cobordism.}
 \vskip 1em

We say that a $G$--manifold $X$ is the $G$--boundary of a $G$--manifold with boundary $W$ if $\partial (W)=X$, and  $X$ has a collared $G$--invariant neighborhood in $W$  of type $X\times [0,1)$ with product $G$--action, which is assumed to be trivial  on the second factor; we will write $X=\partial_G(W)$.
We will also say that $W$ has G--boundary $X$.
Note that, $T^*_G(W)|_{X} = T^*_G(X) \times \Bbb R$.

 The $G$--equivariant operation of restriction 
 to $X$ induces restriction $K_G$--theory homomorphisms
 $$
 \rho^s_{W,X}: K^s_G(T^*_G(W))\to K^s_G (T^*_G(X)\times {\Bbb R}), \quad
 s=0,1.
 $$
 By definition, 
 $$
 K^1_G (T^*_G(X)\times {\Bbb R})
 \cong K^0_G (T^*_G(X)\times {\Bbb R}^2).
 $$
 Hence, for $s=1$,  
 $$
 \rho^1_{W,X}: K^1_G(T^*_G(W))\to K^0_G (T^*_G(X)\times {\Bbb R}^2).
 $$
  
  \proclaim{Definition 2.1 [Ca]}  Let $W$ be a $G$--manifold 
  having as $G$--boundary the $G$--manifold $X$, i.e., $\partial_G(W)=X$. 
  The symbol $G$--boundary map $\partial^W_X: K^1_G(T^*_G(W))\to 
  K^0_G (T^*_G(X))$ is 
  defined by the following equality
  $$
  \partial^W_X = \left( \beta^0_{T^*_G(X)} \right)^{-1}\circ
  \rho^1_{W,X},         
  $$
  where $\beta^0_{T^*_G(X)}: K^0_G (T^*_G(X)) \to K^0_G 
  (T^*_G(X)\times {\Bbb R}^2)$ is the equivariant Bott isomorphism. 
  \endproclaim

  \proclaim{Definition 2.2 [Ca]} Let $\Upsilon_i=(X_i, \sigma_i)$
with $X_i$ a $G$--manifold and $\sigma_i\in K^0_G (T^*_G(X_i))$,
  $i=1,2$.
  Then we say that   $\Upsilon_1$ and $\Upsilon_2$ are symbol
  $G$--cobordant if 
  there exists a pair ${\Cal W} = (W, \sigma)$, with 
$W$ a $G$--manifold with boundary, and 
and $\sigma\in K^1_G (T^*_G(W))$,   
  such that the $G$--boundary of $W$ is
  $X_1 \sqcup X_2$, and
  $$\partial^W_X (\sigma)= - \sigma_1 \oplus \sigma_2.$$
  \endproclaim

  \noindent If $\Upsilon_1$ and $\Upsilon_2$ are symbol $G$--cobordant,
  we will write
  $\Upsilon_1 \sim \Upsilon_2$.  Note that  $\sim $ is an equivalence relation.

  \vskip 1em
  \noindent {\bf 3.  Embeddings, Index, and Symbol $G$--Cobordisms.}
  \vskip 1em

  We will now describe the behaviors of the index, restriction and boundary
  homomorphisms 
  with respect to $G$--embeddings. 
  We will require that all $G$--embeddings
  $\varphi: X\to Y$ of $G$--manifolds 
are $G$--equivariant, $K$--oriented, and 
admit   open $G$--invariant tubular neighborhoods.
  If $X$ and $Y$ are $G$--manifolds with boundary, 
  we also assume that a $G$--embedding
  $\varphi: X\to Y$ restricts to a $G$--embedding on $\partial X$.
(This follows from the existence of tubular neighborhoods for $G$--
manifolds with boundary.)

  Let $\varphi: X\to Y$ be a $G$--embedding of $X$ into $Y$. Then
  $\varphi$ 
  induces $K_G$--theory 'wrong way functoriality'
  shriek maps $\varphi_!^s :K_G^s(T^*_G(X))\to K_G^s(T^*_G(Y))$,
  $s=0,1$,
  defined as below [At].

  \noindent Let $N$ be an open $G$--invariant tubular neighborhood of 
  $\varphi (X)$ in $Y$and let 
  $$
  \tau^s_{X,N}: K_G^s(T^*_G(X))\to K_G^s(T^*_G(N)),\quad s=0,1,
  $$
  be the Thom homomorphism. Moreover, let 
  $$
  \kappa^s_{N,Y}: K_G^s(T^*_G(N))\to K_G^s(T^*_G(Y)),\quad s=0,1,
  $$
  be the $K_G$--theory maps induced by the open 
  embedding $\kappa: T_G^*(N)\to T_G^*(Y)$. 

  \proclaim{Definition 3.1}  
  Let $\varphi:X\to Y$ be a $G$--embedding of $G$--manifolds.
  Then the shriek map 
  $$\varphi_!^s: K_G^s(T^*_G(X))\to K_G^s(T^*_G(Y)), \quad s=0,1,$$
  is defined by the following equality
  $$
  \varphi_!^s=\kappa^s_{N,Y} \circ \tau^s_{X,N}, \quad s=0,1.
  $$
  \endproclaim

  \noindent The following deep theorem, which states 
the invariance of the transverse index under push-forwards, is proved by Atiyah in [At]. 

  \proclaim{Theorem 3.2 [At]}  
  Let $\varphi:X\to Y$ be a $G$--embedding of $G$--manifolds.
  Then the diagram below is commutative.
  $$
  \CD
  K^0_G(T^*_G(X))   @>{\varphi_!^0}>>   K^0_G(T^*_G(Y))\\
  @VV Ind V   @VV Ind V\\
  \Cal D'(G)        @>{Id}>>         \Cal D'(G)   
  \endCD
  $$ 
  \noindent Moreover, if $\varphi: X\to Y$ and $\psi: Y\to Z$ are
  $G$--embeddings of $G$--manifolds, 
  then  
  $$(\varphi \circ \psi)_!^s =\varphi^s_! \circ \psi_!^s, \quad s=0,1.$$
  \endproclaim

  \noindent In the above theorem, $Y$ can also be assumed to be
  non--compact
  because of the functoriality of the shriek maps, and 
  the open $G$--embeddings index invariance, [At].

  The following results, detailing the behavior of the   
  transverse symbol with respect to $G$--embeddings and push-forwards,  are generalizations of analogous results in [Ca]. We omit their proofs.

  \proclaim{Lemma 3.3}  
  Let $\varphi:X\to Y$ be a $G$--embedding of $G$--manifolds.
  Then the diagram below is commutative.
  $$
  \CD
  K^0_G(T^*_G(X))   @>{\varphi_!^0}>>   K^0_G(T^*_G(Y))\\
  @VV{\beta_X^0}V                           @VV{\beta_Y^0}V\\
  K^0_G(T^*_G(X)\times{\Bbb R}^2) @>{{\tilde \varphi}_!^0}>>
  K^0_G(T^*_G(Y)\times \Bbb R^2)    
  \endCD
  $$
  Here 
  $$\beta^0_Z: K^0_G(T^*_G(Z)) \to K^0_G(T^*_G(Z)\times \Bbb R^2), \,
  Z=X,Y, $$
  is the Bott periodicity isomorphism
  and $\tilde {\varphi}$ is the lift of $\varphi$ to $X\times {\Bbb
  R}^2$. (The 
  action of $G$ on ${\Bbb R}^2$ is trivial.)  
  \endproclaim

	  \proclaim{Lemma 3.4}  Let $\varphi:X\to Y$ be a $G$--embedding of
	  $G$--manifolds with boundary. Then the diagram below is
	  commutative,
	  $$
	  \CD
	  K^1_G(T^*_G(X))   @>{\varphi_!^1}>>   K^1_G(T^*_G(Y))\\
	  @VV{\rho_{X, X_0}^1}V               @VV\rho_{Y, Y_0}^1V\\
	  K^0_G(T^*_G(X_0)\times {\Bbb R}^2) @>{\tilde \varphi}_{0\, !}^0>>
	  K^0_G(T^*_G(Y_0)\times {\Bbb R}^2),   
	  \endCD
	  $$
	  where $\partial Z= Z_0$, $Z=X,Y$, $\varphi_0 = \varphi|_{X_0}$,
	  and
	  $\tilde {\varphi_0}$ is the lift 
	  of $\varphi_0$ to $X_0\times {\Bbb R}^2$. (Notation as in
	  Section  2 and Lemma 3.3.)
	  \endproclaim

	  \proclaim{Proposition 3.5}  
	  Let $\varphi:X\to Y$ be a $G$--embedding of $G$--manifolds with
	  boundary.
	  Then the diagram below is commutative,
        $$
	  \CD
	  K^1_G(T^*_G(X))   @>{\varphi_!^1}>>   K^1_G(T^*_G(Y))\\
	  @VV\partial^X_{X_0}V
	  @VV\partial^Y_{Y_0}V\\
	  K^0_G(T^*_G(X_0))        @>{ \varphi_{0\, !}^0}>>
	  K^0_G(T^*_G(Y_0))    
	  \endCD
	  $$
	  where $\varphi_{0}={\varphi|_{X_0}}$, $\partial Z= Z_0$, 
	  and the boundary map $\partial^Z_{Z_0}$, 
	  is as in  Definition 2.1. ($Z=X,Y$.)
	  \endproclaim

        \vskip 1em
	  \noindent {\bf 4.  The Main Result.}
	  \vskip 1em

In this final section,  we will prove our main results,
Theorem 4.3 and Theorem 4.6.
 But first two lemmas.

 Let $M$ be a $G$--manifold and $\sigma 
	  \in K_G(T^*_G(M))$. Recall that $\Upsilon=(M,\sigma)\sim 0$
(that is, $\Upsilon$ is  symbol $G$--cobordant to zero)
	  if there exist a $G$--manifold with boundary $W$ with $X=\partial_G(W)$,
	  and 
$\omega \in  K_G^0(T^*_G(W)$ such that
	  $\partial^W_X(\omega)=\sigma$.

	  \proclaim{Lemma 4.1}  Let $M$ be a $G$--manifold and  $\sigma 
	  \in K_G(T^*_G(M)$, with $\Upsilon=(M,\sigma)\sim 0$ via 
the $G$--cobordism $W$. 
	  If $\varphi : W\to Y$ is a $G$--embedding of $W$ into the 
	  $G$--manifold with boundary $Y$, with restriction to the 
	  boundary given by $\varphi_0$, then 
	  $$\partial^Y_{\partial Y}(\varphi^1_{!}(\omega))= \varphi_{0\,
	  !}^0(\sigma),$$ 
	  which implies $\Upsilon_Y=(\partial Y,
	  \varphi_{0\,!}^0(\sigma))\sim 0$.
	  \endproclaim

	  \noindent {\it Proof.} Apply the results of Section 3. $\qed$

\proclaim{Lemma 4.2} We have, 
$$
K^i_G(T^*_G(E))=0,\ i=0,1,
$$
where $E$ is equal to the $G$--module $[0,1) \times {\Bbb R}^t$, $t>0$. 
\endproclaim

\noindent {\it Proof.}  Firstly, $K^i_G(T^*(E))=0$, $i=0,1$,
since this algebra is 
$G$--contractible, [Ca].    
Next, $K^i_G(T^*_G(E))$, $i=0,1$, 
can be decomposed as the direct sum 
$K_G^i(T^*_G(E))= \oplus_j K_G^i(T_G^*(E(j)-E(j+1)))$,
where $E(j)=\{ x\in E : \hbox{dim } G_x \geq j\}$, $i=0,1$, 
see [At], Theorem 8.4. 
Also note that  that a decomposition 
similar to the above can be proved for 
$K_G^i(T^*(E))$, $i=0,1$, with similar methods.
Of course, in this latter case,
each of the factors is 
the zero module.

Since on each of the  spaces appearing in the above decompositions
the action of the group can be assumed to have only finite stabilizers, 
by the Bott periodicity theorem
it follows that each of the factors in the
two decompositions are isomorphic in pairs,  from which 
the result follows.  
Note that the Bott periodicity theorem can be 
applied to this case where the action is non-trivial by
[Ph], Remark 2.8.7. $\qed$

	  \proclaim{Theorem 4.3} Suppose that  $M$ is a $G$--manifold, and let $D$  be a $G$--p.d.t.e. operator on $M$ with $G$--transverse symbol class 
	  $\sigma\in K_G^0(T^*_G(M))$. If  
	  $\Upsilon=(M, \sigma)\sim 0$, then 
	  $Ind(D)=0$.
	  \endproclaim

	  \noindent {\it Proof.} Let ${\Cal W}=(W, \omega)$, for some 
	  $\omega\in K^1_G(T^*_G(W))$, 
	  be a symbol $G$--cobordism between $\Upsilon$ and 0.
	  Let $\varphi: W\to E$ be a $G$--embedding of $W$ into 
	  $E=[0,1) \times {\Bbb R}^t$ (see  e.g. [Ca]).
	  If  we denote by $\varphi_0$
the restriction of $\varphi$ 
to the boundary, then by Lemma 4.1, 
	  $$\partial^F_{\partial F}(\varphi^1_{!}(\omega))= \varphi_{0\,
	  !}^0(\sigma). $$
	  By Lemma 4.2 and the results in Section 3, 
	  $\varphi_!^1(\omega)=0$, which implies $Ind(D)=0$. $\qed$

	  \noindent As a corollary, we have,

	  \proclaim{Corollary 4.4}  Suppose that  $M_i$ is a 
	  $G$--manifold, and let $D_i$ 
	  be a $G$--p.d.t.e. on $M_i$ with $G$--transverse symbol class 
	  $\sigma_i\in K_G^0(T^*_G(M_i))$, $i=1,2$. Assume that 
	  $\Upsilon_1=(M_1, \sigma_1)\sim \Upsilon_2=(M_2, \sigma_2)$.
	  Then 
	  $Ind(D_1)=Ind(D_2)$.
	  \endproclaim

	  \noindent {\it Proof.} We have that $$(M_1\sqcup M_2,
	   -\sigma_1\oplus \sigma_2)\sim 0.$$ Now 
	   $(-\sigma_1\oplus \sigma_2)$ is the $G$--transverse symbol of
	   $D_1^*\oplus D_2$. Then the claim follows from Theorem 4.3.
	   $\qed$
	   \vskip 1em

We can now prove the invariance under cobordism of the orbifold index.

\proclaim{Definition 4.5} Let $Q_i$, be an orbifold. Assume that 
$Q_i$ arises as the locally free quotient of the $G$--manifold
$M_i$, $i=1,2$. Let $P_i$ be an elliptic pseudo-differential 
operator on $Q_i$ with symbol $p_i$, and let  ${\tilde P}_i$  be its lift to $M_i$ with symbol $\tilde{p}_i $, $i=1,2$.
Then  ${\tilde P}_i$ is a $G$--p.d.t.e. on $M_i$ with $G$--transverse symbol class $\sigma_{\tilde{P}_i}\in K_G^0(T^*_G(M_i))$, $i=1,2$.
We say that $(Q_1, p_1)$ is orbifold symbol cobordant to $(Q_2, p_2)$
if $\Upsilon_1=(M_1, \sigma_{\tilde{P}_1})\sim \Upsilon_2=(M_2, \sigma_{\tilde{P}_2})$.

\endproclaim

\proclaim{Theorem  4.6} Let $Q_i$, be an orbifold and   let $P_i$ be an elliptic pseudo-differential operator on $Q_i$ with symbol $p_i$, 
$i=1,2$. Suppose that $(Q_1, p_1)$ is orbifold 
symbol cobordant to $(Q_2, p_2)$. Then  $Ind(P_1)=Ind (P_2)$.
\endproclaim

We would like to thank the Department of Mathematics of the
University of Florence, Italy, for their
warm hospitality during the period while this note was written, and the 
University of Colorado Sabbatical
Program. We would also like to thank an anonymous referee for 
detecting an essential gap in our original proof and
for detailed suggestions that greatly improved the exposition.

\vskip 1em
 \noindent {\bf References}
\itemitem {[At]} M.F. Atiyah, Elliptic operators and compact
     groups, Lecture Notes in Mathematics no. 401, Springer, Berlin, 1974.
	   \itemitem {[Bl]} B. Blackadar,  $K$-theory for 
	   operator algebras, second edition, Mathematical Sciences 
	   Research Institute Publications, no. 5, Cambridge University Press, 
	   Cambridge, 1998. 
	   \itemitem {[Br1]} M. Braveman, New proof of the cobordism
	   invariance of the index, 
	   Proc. Amer. Math. Soc. 130 (2002), 1095--1101.
	    \itemitem {[Br2]}   M.  Braveman, Cobordism invariance of the 
		index by a transversally elliptic operator, Appendix J, in 
V. Guillemin, V. Ginzburg and Y. Karshon,
Momentum maps, cobordisms, and Hamiltonian group actions,
Mathematical Surveys and Monographs, vol. 98, American Mathematical
Society, Providence, 2002.
	   \itemitem {[BGV]} N. Berline, E. Getzler and M. Vergne, Heat
	   kernels and
	   Dirac operators, Grundleheren der Mathematical Wissenshaften no. 298,
	   Springer--Verlag, Berlin 1992.
	   \itemitem {[BV]} N. Berline and M. Vergne, The Chern character of
	   a  transversally elliptic symbol and the equivariant index, Inv.
	   Math. 124 (1996), 11--49.
	   \itemitem {[Ca]} C. Carvalho, A $K$--theory proof of the
	   cobordism 
	   invariance of the index, 
	   $K$-Theory 36 (2005), 1--31 (2006). 
	   \itemitem {[Dr]} K.S. Druschel, Oriented orbifold cobordism,
	   Pacific J.
	   Math. 164 (1994), 299--319.
	   \itemitem {[Fa1]} C. Farsi, $K$-theoretical index theorems for
	   orbifolds,
	   Quat. J. Math. 43 (1992), 183--200.
	   \itemitem {[Fa2]} $\underline {\hskip 1.0in}$   Orbifold spectral
	   theory,
	   Rocky Mtn. J.  31 (2001).
	   \itemitem {[Fa3]} $\underline {\hskip 1.0in}$ Orbifold
	   $\eta$-invariants, 
	   Indiana Univ. Math. J. 56 (2007), 501--521.
\itemitem {[Hi]} N. Higson, The signature theorem for
	   $V$-manifolds,
	   Topology 30 (1991), 439--443.
	   \itemitem {[Mo]} S. Moroianu, Cusp geometry and the cobordism 
	   invariance of the index,
	   Adv. Math. 194 (2005),  504--519.
	   \itemitem {[Kw1]} T. Kawasaki,  The signature theorem for
	   $V$--manifolds,
	   Topology 17 (1978), 75--83.
	   \itemitem {[Kw2]} $\underline {\hskip 1.0in}$   The Riemann Roch
	   Theorem
	   for complex $V$-manifolds, Osaka J. Math. 16 (1979), 151--159.
	   \itemitem {[Kw3]} $\underline {\hskip 1.0in}$   The index of
	   elliptic
	   operators over $V$-manifolds, Nagoya Math. J. 84 (1981), 135--157.
	   \itemitem {[Pa]} R. Palais, Seminar on the Atiyah-Singer index
	   theorem,
	   Annals of Mathematics Studies, No. 57, Princeton University
	   Press, 
	   Princeton, 1965.
	   \itemitem {[Ph]} N. C. Phillips,  Equivariant K-theory
	   and freeness of group actions on $C^*$--algebras,
Lecture Notes in Mathematics, no. 1274, Springer-Verlag,
Berlin   and Heidelberg,  1987. 
   \itemitem {[Ve]} M. Vergne, Equivariant index formulas for
	   orbifolds,
	   Duke Math.  J. 82 (1996), 637--652.\end